# Damping Tuning Considering Random Disturbances Adopting Distributionally Robust Optimization

Yuhong Wang, *Senior Member, IEEE*, Xinyao Wang, *Student Member, IEEE*, Chen Shen, *Senior Member, IEEE*, Jianquan Liao, *Member, IEEE*, Qianni Cao, *Student Member, IEEE,* Yufei Teng, *Member, IEEE,* Huabo Shi, Gang Chen

*Abstract*—In scenarios where high penetration of renewable energy sources (RES) is connected to the grid over long distances, the output of RES exhibits significant fluctuations, making it difficult to accurately characterize. The intermittency and uncertainty of these fluctuations pose challenges to the stability of the power system. This paper proposes a distributionally robust damping optimization control framework (DRDOC) to address the uncertainty in the true distribution of random disturbances caused by RES. First, the installation location of damping controllers and key control parameters are determined through Sobol sensitivity indices and participation factors. Next, a nonlinear relationship between damping and random disturbances is established with Polynomial Chaos Expansion (PCE). The uncertainty in the distribution of disturbances is captured by ambiguity sets. The DRDOC is formulated as a convex optimization problem, which is further simplified for efficient computation. Finally, the optimal control parameters are derived through convex optimization techniques. Simulation results demonstrate the effectiveness and distribution robustness of the proposed DRDOC.

## I. Introduction

As the installed capacity of renewable energy sources (RES) has grown rapidly worldwide, their penetration into power grids has also increased significantly. Advanced transmission technologies, such as High-voltage direct current transmission based on voltage source converter (VSC-HVDC), are regarded as one of the ideal solutions for enabling the long-distance integration of RES [1]. However, the intermittency and increasing uncertainty of RES, along with the control dynamics of VSC-HVDC, present new challenges to the stability and reliability of power system [2].

The steady-state operating point of power system with significant integration of RES, such as wind power, typically exhibits random variations [3]. Traditional damping controllers, which rely on local information about the operating point, provide a certain degree of robustness around this point. This approach has two key limitations. First, the effectiveness of a damping controller is influenced by the choice of installation location and the selection of input and output signals. Second, these controllers may struggle to handle the large fluctuations in wind power generation [4]-[5]. In weak power grids, the pronounced variability of wind power generation can cause drastic shifts in the operating point, increasing the risk of system instability.

The installation of damping controllers and the selection of control signals have not been systematically explored in power grids with a high proportion of RES. In AC grids, the installation locations of Power System Stabilizers (PSS) are typically determined by calculating the residual coefficients of synchronous generators [6]. However, this method may not be directly applicable to modern grids with substantial RES integration. Some studies suggest installing damping controllers at the point of common coupling (PCC), but this approach is limited to point-to-point transmission scenarios. For signal selection, Ref. [7] propose using the voltage or current at the PCC, though it lacks a solid theoretical foundation. Additionally, the choice of control parameters also impacts system stability. Ref. [8] applies Particle Swarm Optimization (PSO) to optimize all PSS control parameters. However, this method may increase controller losses and potentially lead to instability in other modes. Therefore, it is essential to develop methods that can theoretically determine the most effective installation locations, signals, and parameters for damping controllers.

Extensive research has been conducted to address power system instability caused by random disturbances. Two commonly adopted control methods are robust optimization (RO) [9] and stochastic optimization (SO) [10]. RO determines the optimal decision by considering the worst-case scenario, where uncertain parameters are set within specified ranges. However, this approach tends to be overly conservative, as it does not account for the unknown probability distribution of the uncertain parameters [11]. SO assumes that the probability distribution of uncertain parameters is known and derives the optimal decision accordingly. While this method avoids the conservatism of RO, the probability distribution of wind power volatility in power systems is often unknown. Although some researchers

This work was supported by Smart Grid-National Science and Technology Major Project (2024ZD0801600). *(Corresponding author: C.Shen)*

X. Wang, Y. Wang and J. Liao(corresponding author) are with the College of Electrical Engineering, Sichuan University, Chengdu, China (e-mail: xinyao@stu.scu.edu.cn; yuhongwang@scu.edu.cn; jquanliao@scu.edu.cn).

C. Shen and Q. Cao are with the Electrical Engineering Department, Tsinghua University, Beijing, China (e-mail: shenchen@mail.tsinghua.edu.cn; cqn20@mails.tsinghua.edu.cn).

Y. Teng, H. Shi and G. Chen are with the State Grid Sichuan Electric Power Research Institute, Chengdu, China (e-mail: yfteng2011@163.com; gangchen_thu@163.com; shbo87@163.com).

attempt to estimate this distribution through forecasting and other techniques, discrepancies between predicted and actual values persist and are difficult to quantify [12]. As a result, there is growing interest in distributionally robust optimization (DRO) that combines the strengths of both RO and SO. DRO enhances the robustness of optimization decisions under uncertainty by considering ambiguity sets of possible uncertainty distributions, rather than relying on a precise distribution model [13]. Ref. [14] discusses the principle of DRO, the construction of ambiguity sets, and how to handle reconstruction.

DRO is typically formulated as a nonlinear, non-convex optimization problem, making it complex to solve and challenging to achieve the global optimal solution. Ref. [15] transforms distributionally robust control (DRC) into a semi-definite program (SDP) through the reconstruction of specific chance constraints. Ref. [16] employs conditional value at risk to derive the boundary. However, these solution methods are tailored to specific optimization problems or scenarios and lack a universal strategy applicable to all DRO formulations. For damping optimization control, the objective function generally focuses on improving the damping of key modes. This necessitates an accurate characterization of the complex nonlinear relationship between damping and disturbances [17]. Ref. [18] approximates the damping distribution using second-order sensitivity, but this method struggles to capture variation characteristics in highly complex nonlinear relationships. Ref. [19] adopts a two-point estimation method to approximate the damping distribution. However, its accuracy is limited when dealing with nonlinear functions. In recent years, polynomial chaos expansion (PCE) has been applied to damping evaluation [20]. PCE effectively captures complex nonlinear relationships while requiring fewer data points. Ref. [21] applies PCE to assess the impact of illumination and load on damping. Furthermore, the mathematical properties of PCE allow for straightforward computation of the expected value and variance of damping, which is advantageous for solving distributionally robust damping optimization problems.

To address the issues and limitations of existing methods, this paper proposes a damping tuning strategy adopting DRO. It aims to design additional damping controllers to improve the system damping from the uncertainty of RES. First, a method for determining the installation location of a damping controller, input and output signals, and key control parameters is proposed based on Sobol sensitivity indices and participation factors. Second, we consider the deviation between the actual injected power and the reference power of RES as disturbance. the nonlinear relationship between damping and random disturbances is established based on PCE. Furthermore, the uncertainty in the true distribution of the disturbance is described using ambiguity sets. The uncertainty damping optimization problem is modeled as a distributionally robust optimization problem and transformed into a solvable form. Finally, the optimized control parameters are obtained by convex optimization solving algorithm.

The main contributions of this paper are as follows. 1) To improve the uncertainty damping affected by RES, a distributionally robust damping optimization control framework (DRDOC) is proposed. This framework employs PCE to model the nonlinear relationship between damping and random disturbances and integrates it into the DRO problem. The non-convex optimization problem with ambiguity set constraints is effectively solved. 2) The analytical process of simplifying the DRO into a solvable form is derived, utilizing the mathematical properties of PCE and DRO. This reduces the non-convex optimization problem to a convex one, significantly improving computational feasibility. 3) A method for selecting the installation locations, signals, and parameters of damping controllers is proposed, based on Sobol sensitivity analysis. This method is scalable to cases involving multiple parameters.

The remaining sections of this paper are organized as follows. Section II presents the method for selecting the installation location and control parameters of the damping controller. Section III develops a polynomial approximation of damping, considering random disturbances based on PCE. Section IV provides the analytical form and simplification process of the DRDOC, along with its overall framework. Section V validates the effectiveness and distributional robustness of the DRDOC on a wind power system connected to AC grids via multi-terminal flexible DC based on VSC (VSC-MTDC) in MATLAB/Simulink. Section VI concludes the paper.

## II. INSTALLATION LOCATION AND SENSITIVE PARAMETERS

This section introduces the structure of the additional damping control. Eigenvalue sensitivity analysis is applied to determine the installation location and input signals. The Sobol method is employed to select parameters within the stabilizer that are sensitive to weakly damping modes. These sensitive parameters are then optimized within the DRDOC framework to maintain system stability under random disturbances.

### A. Control Structure of the additional damping control

PSSs are generally installed on generators to suppress oscillations in AC grids. They enhance the weak damping characteristics of system near the oscillation frequency by adding an additional control loop, without impacting the steady-state operation of the power system. Based on this idea, the control structure of the additional damping control is illustrated in Fig. 1. The gain $K_m$ and time constants $T_1$, $T_2$, $T_3$, $T_4$ are the parameters to be tuned.

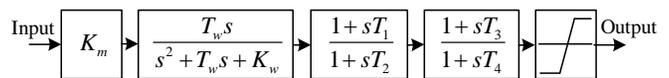

Fig. 1 The control structure of damping control.

It can be synthesized to a single state space representation as follows:

$$\qquad(1)$$

where $x_c$ are the state variables of the damping control. $y_c$ and $u$ are the output and input vectors of the damping control. The matrices $A_c$, $B_c$, $C_c$, $D_c$ are exclusively calculated from the damping controllers' parameters so that they are parameterized by the tunable parameter vector , as follows:

$$\qquad(2)$$

All parameters in $G_{PSS}$ are contained in . Based on Eq. (2), $A_c$, $B_c$, $C_c$, $D_c$ have no relationship with the power system operation.

*B. Influence of Installation Location.*

In traditional AC grids, PSSs are installed on synchronous generators. The angular frequency deviation or active power deviation is processed through a transfer function and fed back to the excitation circuit as an additional excitation signal. However, in MTDC grids, the oscillation mechanism is different from that of the synchronous generators. Therefore, the installation location and excitation signal of the controller cannot be determined based on practices in AC grids. Based on participation factors, this paper outlines the process for selecting the most critical locations and excitation signals associated with oscillation modes.

The state space of the MTDC system can be expressed as:

$$\qquad(3)$$

where $A$, $B$ and $C$ are state matrix, input matrix and the output matrix. Eigen analysis of matrix $A$ produce system eigenvalues $\lambda_i$. The location is determined by using participation factors (PF).

PF reflects the extent to which a specific state contributes to a particular mode. The involvement of state $s$ in the $i$th mode can be expressed as:

$$\qquad(4)$$

where $t_{is}$ and $v_{is}$ are right and left eigenvectors.

Hence, converter with the highest participation factors in the targeted mode are ideal candidates for implementing control through damping controller.

The state space of the MTDC system equipped with additional $l$ damping controller can be expressed as:

$$\qquad(5)$$

where $B_j$ and $C_j$ are the column-vector input matrix and the row-vector output matrix corresponding to the $j$th controller, respectively.

The open loop transfer function corresponding to a specific input/output is:

$$\qquad(6)$$

Where is the residue with respect to the $i$th mode and it is given by

$$\qquad(7)$$

When a damping controller is installed in the $j$th converter, the $i$th eigenvalue variation is calculated as follows:

$$\qquad(8)$$

where is phase lead-lag compensators of the controller.

The location with the largest residue for the oscillation mode of interest is an optimal candidate for installing additional control. Eq. (8) provides a quantitative measure of the influence of the additional control on the targeted mode and other modes. It would be utilized as a design constraint in Section IV. The DRDOC aims to minimize its impact on other modes, ensuring that additional control does not adversely affect the stability characteristics of those modes.

*C. Sensitivity Analysis of Control Parameters.*

The PF-based parameter selection process evaluates the effect of state variables on modal damping. However, it cannot directly determine the relationship between parameters and damping. To address this limitation, the Sobol method is used to efficiently extract the most sensitive parameters. To address this limitation, the Sobol method is employed to efficiently extract the most sensitive parameters.

The Sobol global sensitivity analysis based on variance is a widely used method for identifying key factors. Its core idea is to decompose the target model into functions of parameter combinations. This allows the impact of individual parameters or parameter sets on the power system response to be evaluated through their variance. The principle is explained as follows.

The variance of system response output can be expressed as a function of a series of inputs as follows:

$$\qquad(9)$$

with

$$\qquad(10)$$

where $v_s$ and $z_s$ are output and input variables of system. $z_{s\sim i}$ is the vector of all variables except $z_{si}$. $V_i$ and $V_{ij}$ are the variance of the corresponding conditional expected value.

The variance expression of the model output allows the influence of input variables on the output $v_s$ to be represented using the Sobol index. Sobol indices, derived from variance decomposition, are widely used to measure the sensitivity of input variables [22]. The total index captures the overall importance of each input, including its direct effects and higher-order interactions. According to these definitions, input variables with larger Sobol indices are more sensitive.

The total index $S_{Ti}$ is:

$$\qquad(11)$$

Based on Eqs. (9-11), the influence of all control parameters on the small signal stability of the system can be quantified. Hence, the most sensitive parameters are selected, enabling effective optimization of damping control.

### III. DAMPING PROBABILITY DISTRIBUTION CONSIDERING RANDOM DISTURBANCES

When large-scale wind farms are connected to grids, the inherent randomness of wind power generation may cause the system's equilibrium points to shift. In a weak power grid, it may lead to instability under disturbances. For a specific system with fixed operating mode, topology, nodes and control parameters, the power flow is uniquely determined. Consequently, the equilibrium points and eigenvalues are also determined.

From this perspective, the eigenvalues and damping could be regarded as functions of the node injected power and control parameters. However, due to the complexity of power system equations, this relationship is inherently nonlinear. When the node injected power and control parameters are treated as random variables, the damping becomes functions of the random variable. Thus, this section constructs probability density function of the system damping, taking the randomness into account.

*A. Polynomial Approximation of Critical Damping.*

To address the weak damping issue caused by random disturbances, this section focuses on rapidly calculating the probability distribution of system damping with sufficient accuracy.

The system can be written in the form of a general differential-algebraic equation system:

$$ \quad (12)$$

where and are the differential and algebraic equations respectively, and are the state variables and algebraic variables of the system. includes the state variables of wind turbine and MTDC. are the parameters of dimension $d$, encompassing the operating parameters and control parameters that influence the operating point. In this paper, is the error between the actual power injected by the node and the injected power reference value.

For any operating point , the small signal stability of the system can be determined by the linearized model near the operating point:

$$ \quad (13)$$

where is the increment of system state variables, and is the system state matrix at the operating point.

The randomness of wind speed can cause the operating point to shift, potentially leading to instability. The relationship between system damping and parameters is a complex implicit function, making its analytical expression difficult to derive. However, this relationship can be approximated using a polynomial function.

The polynomial approximation is expressed as a linear combination of polynomial basis functions:

$$ \quad (14)$$

where represents the damping calculation equation, is its chaotic polynomial approximation, is the polynomial basis, is the coefficient of the kth polynomial basis, and is the number of polynomial basis.

The polynomial basis are generally selected as the product of various univariate polynomials as follows:

$$ \quad (15)$$

where represents polynomial of degree $n_i$ with respect to the single variable $p_i$. In the theory of PCE, the commonly used polynomial basis belongs to the orthogonal polynomial system within the Askey polynomial. The weight functions corresponding to these polynomials are the PDF of common distributions, after appropriate normalization and translation. This establishes a connection between the orthogonal polynomial system and the probability distribution.

The multidimensional Hermite PCE representing the normal distribution in Askey can be expanded as follows:

$$ \quad (16)$$

where is the independent standard normal random variables, which can be obtained by Nataf transformation. $n$ represents the degree of freedom of the random input uncertainty, which refers to the minimum number of variables required to describe the random input. is a Hermite polynomial with $m$ degrees of freedom:

$$ \quad (17)$$

By employing a finite $m$-order Hermite chaos polynomial expansion, we obtain unknown $N_a$ constant coefficients and corresponding equations, $N_a$ can be obtained as:

$$ \quad (18)$$

The truncated PCE is expressed as:

$$ \quad (19)$$

By applying the Monte Carlo method to the proxy model constructed with PCE, the PDF of the damping ratio can be efficiently obtained.

*B. PCE Error.*

The approximation error between any polynomial function and the original function is defined as the weighted $L^2$ norm error of the difference between the two:

$$ \quad (20)$$

The error of the PCE model primarily arises from two sources. The first is the truncation error, which occurs when the model is truncated at the $m$-th order, excluding higher-order polynomial terms. The second is the coefficient estimation error. is the best polynomial approximation of the orthogonal projection, then the error can be written as:

$$\text{(21)}$$

Based on Eq. (20-21), the polynomial basis in this paper is the Hermite basis, therefore and are equivalent.

### IV. DRDOC FRAMEWORK

An ambiguity set is constructed for distribution of disturbance, with controller parameters defined as decision variables, forming the basis of the proposed DRDOC framework.

*A. Ambiguity Sets Based on Wasserstein Distance.*

Although the true probability of random variables is unknown, it is reasonable to assume that it closely approximates the empirical distribution derived from historical data fitting. Based on this assumption, a set of probability distributions is constructed around the empirical distribution to ensure that the true distribution is captured within this set. To measure the distance between these uncertain probability distributions, the Wasserstein distance is adopted.

Distributed robust optimization based on the Wasserstein distance (W-DRO) offers two significant advantages. First, compared to other ambiguity sets, it demonstrates superior tail performance and generalization capability. Second, the uncertainty set defined by the Wasserstein distance can include the true distribution of the data with high confidence.

The Wasserstein distance is defined as follows:

$$\text{(22)}$$

where is the set of joint probability distributions of with a marginal distribution of and , is a sample from the set of random variables, is a distance function in the metric space (such as Euclidean distance), and represents the set of random variables. represents the first-order Wasserstein distance between the probability distributions and . It can effectively measure the difference between the two distributions.

Based on historical data and simulation data, an empirical distribution $P_0$ can be given. Wasserstein distance is adopted to define the ambiguity sets:

$$\text{(23)}$$

$P_0$ is an empirical distribution that can either be derived from fitting historical data or represent an unknown distribution characterized by key properties such as mean and covariance, which describe the disturbance. The error between the actual power injected at the node and the reference power is assumed to follow a normal distribution defined by its mean and variance [9]. $\delta$ represents the threshold of the Wasserstein distance, defining the radius of the ambiguity set. The value of $\delta$ determines the level of conservatism in the control decision.

To ensure that the true distribution is within the ambiguity set, the confidence level is chosen to be 1-$\beta$. $\delta$ can be calculated as:

$$\text{(24)}$$

where $N$ is the number of historical samples and $D$ is a constant that can be obtained by solving the following optimization problem:

$$\text{(25)}$$

where $\mu$ is the sample mean.

Ref. [14] demonstrated that as the sample size approaches infinity, the worst-case distribution converges to the true data distribution. Consequently, the optimal decision under the worst-case distribution aligns with the true optimal parameters of the model.

*B. Proposed Optimization Method.*

DRO is a generalization of robust optimization, requiring only that the distribution of uncertain parameters satisfies certain constraints. The uncertainty set is defined as a collection of probability distributions that conform to specific known information. To ensure the distributional robustness of the additional damping controller, the ambiguity set is constructed as shown in Eq. (23). The design problem for the controller, aimed at suppressing negative damping oscillations, is then transformed into a parameter optimization problem under the worst-case damping condition induced by the ambiguity set.

The damping optimization problem under disturbance can be formulated as the following optimization problem:

$$\text{(26)}$$

$$\text{(27)}$$

$$\text{(28)}$$

$$\text{(29)}$$

where and are the lower and upper bounds of , respectively; represents the change of all eigenvalues modes except the weakly damping mode.

Eq. (26) focuses on maximizing the expected worst-case damping under the influence of the ambiguity set. This is accomplished by optimizing the control parameter design. The formulated optimization problem is highly relevant to damping control design, as the controller effectively regulates the critical mode to ensure satisfactory damping under uncertain random disturbance distributions, thereby achieving robustness against variations in disturbance patterns. Additionally, the impact of the additional control on other modes is constrained to remain minimal, preventing adverse effects on these modes.

By leveraging the polynomial fitting properties of PCE, Eq. (26) can be transformed into a simpler form that is easier to solve. The expectation in PCE is computed by integrating the probability distribution of the output random variable. Due to

the linearity of PCE, the expected value can be calculated as:

$$\qquad(30)$$

Due to the orthogonality of the polynomials, the inner product of the constant term and other polynomial basis equals zero:

$$\qquad(31)$$

In view of Eq. (31), Eq. (32) can be rewritten as follows:

$$\qquad(32)$$

The expectation of any polynomial, apart from the constant term, is zero. As a result, the expected value depends entirely on the coefficient of the constant term. Consequently, Eqs. (29-32) can be transformed into:

$$\qquad(33)$$

$$\qquad(34)$$

$$\qquad(35)$$

By applying linear relaxation, the 1-Wasserstein distance is reformulated as a convex constraint [14]. Consequently, the DRDOC is transformed into a convex optimization problem, which can be efficiently solved using standard software packages such as Gurobi [23].

*C. Flowchart of DRDOC Framework.*

Fig. 2 illustrates the overall schematic of the DRDOC framework. Specifically, the main control parameters are initially selected through Sobol sensitivity analysis, while installation locations and signals are determined via factor analysis. Next, data samples are generated by randomly combining disturbances and parameters.

Subsequently, PCE is utilized to model the relationship between disturbances and the damping of the dominant mode, enabling the estimation of the damping PDF with a limited number of samples. The uncertainty in the distribution of disturbance is captured by ambiguity sets. The DRO problem is formulated and further simplified into a solvable form. Finally, the optimized parameters and damping tuning results are obtained by solving the simplified optimization problem.

V. CASE STUDIES

This section evaluates the performance of the DRDOC using a VSC-MTDC model developed in MATLAB/Simulink. The model topology is depicted in Fig. 3, with the parameter settings provided in Table I in Appendix A. First, the locations and parameters are determined through the PF method and Sobol sensitivity analysis. Next, the accuracy and effectiveness of the PCE are validated. Subsequently, the performance of the DRDOC is compared with traditional damping control to demonstrate its effectiveness. Finally, the distributional robustness of the DRDOC is assessed by comparing it with RO and SO.

Fig. 3 Topology of MTDC.

*A. PF and Sobol analysis.*

The weakly damped mode of the system is identified based on the eigenmatrix. The eigenvalue of the weakly damping mode is -2.14 ± i137.84, with a damping ratio of 0.015. The participation factors of this mode in the small-signal model linearization of the VSC-MTDC are shown in Fig. 4. Among the variables, vector control $Z_{P3}$ and DC voltage $V_{dc3}$ of $VSC_3$ have the strongest participation in the weakly damping mode. $I_{1\_L2}$ and $I_{2\_L2}$ represent the transmission dynamics, which have minimal influence on the weak damping mode.

Fig. 4 Participation Factor of weak damping mode.

Fig. 2 Flowchart of DRDOC framework.

Therefore, the DRDOC is implemented at the $VSC_3$ converter station. The variations in DC voltage and active power are selected as the feedback and modulation signals for the damping controller. The control structure in MTDC is shown in Fig. 5.

Fig. 5 Control structure in MTDC.

Subsequently, the Sobol total indexes of control parameters in vector $\boldsymbol{p}_c$ are calculated to determine the parameter that most affects weak damping. The value ranges of the parameters are provided in Table I. To account for differences in the magnitudes of the input parameters, their impact on the output is minimized by normalizing the parameter ranges. The Sobol total indices, derived from 100 sampling points, are presented in Fig. 6. Among the parameters, $K_m$ and $T_1$ exhibit higher Sobol indices compared to others, indicating their greater influence on weak damping. Therefore, $K_m$ and $T_1$ are selected as decision variables in the subsequent steps to enhance the system's robustness against the uncertainty in the true distribution of disturbances.

TABLE I
PARAMETER VALUE RANGE

| Parameters | Value Range |
|---|---|
|  | 1.0-30.0 |
|  | 0.1-1.0 |
|  | 0.01-0.1 |
|  | 0.1-1.0 |
|  | 0.01-0.1 |

Fig. 6 Sobol Total Index of 100 samples.

*B. PCE accuracy.*

The nonlinear relationship between damping and random disturbances is modeled using PCE, which quantifies the probability distribution of damping under the influence of random disturbances. The probability density functions (PDFs) obtained from PCE with a truncation order of 4 and Monte Carlo (MC) simulations are shown in Fig. 7. The system damping probability distribution ranges from −0.02 to 0.08, with the highest probability occurring between −0.02 and 0.02. Generally, damping values below 0.05 are classified as weak damping, while values below 0 are considered negative damping. This clearly indicates that the system is prone to instability due to the influence of disturbances.

The accuracy of the PCE is further validated using the root mean squared error (RMSE) and average absolute error (AAE). The calculated coefficients are as follows:

$$ \tag{36}$$

$$ \tag{37}$$

The errors of PCE with different truncation orders are presented in Table II. As the order of the PCE increases, its accuracy improves. However, this improvement comes at the cost of increased simulation time. Despite this, PCE can capture high-order nonlinear dynamics while offering significantly faster computational speed compared to MC simulations.

Fig. 7 Damping distribution of PCE.

TABLE II
ERROR EVALUATION

| METHOD | $\varepsilon_{RMSE}$ | $\varepsilon_{AAE}$ | Time (s) |
|---|---|---|---|
| MC | / | / | 210.30 |
| 2 order PCE | 0.1152 | 0.1749 | 5.66 |
| 3 order PCE | 0.0671 | 0.0857 | 7.64 |
| 4 order PCE | 0.0091 | 0.0176 | 10.40 |
| 5 order PCE | 0.0057 | 0.0116 | 17.27 |
| 6 order PCE | 0.0031 | 0.0092 | 44.62 |

*C. Effectiveness of DRDOC.*

After identifying the key control parameters and obtaining the damping probability distribution, the DRO problem can be solved. The real distribution of uncertain disturbances is characterized by the Wasserstein distance. Historical data is sourced from a wind farm (ID: 4209) in the National Renewable Energy Laboratory (NREL) database [24]. Fig. 8 depicts the empirical distribution $P_0$, constructed using data with a reference power of 0.5 p.u..

The worst-case disturbance distribution, which results in the most detrimental damping, is selected from within the ambiguity sets, and the optimized parameters are presented in Table III. In this paper, we design the damping optimization for random disturbances in a single wind farm. For scenarios where multiple wind farms serve as independent random disturbances within the system, the proposed method remains broadly applicable.

Fig. 8 Empirical distribution of power prediction error.

The eigenvalue distribution of the system, both before and after the addition of the damping controller, is shown in Fig. 9. The weakly damping mode shifts further into the left half-plane, becoming more stable. Meanwhile, changes in other modes are minimal or negligible, indicating that the proposed method does not cause mode drift or negatively affect system stability.

TABLE III
OPTIMIZED PARAMETER

| Parameters | Initial value | Optimized value |
|---|---|---|
|  | 1.0-30.0 | 9.6 |
|  | 0.1-1.0 | 0.415 |

Fig. 9 Characteristic roots of the initial system and system with DRDOC.

Furthermore, the effectiveness of the additional damping controller designed using the traditional method is compared with that of the DRDOC. Ref. [25] proposed a damping controller design method suitable for MTDC, featuring a basic structure composed of a band-pass filter and a phase compensation stage connected in series. The transfer function of the traditional PSS-type damping controller, designed following the method in [25], is expressed as:

$$ \tag{38} $$

The traditional damping controller is placed at the same location as DRDOC, with identical input and output signals. In the simulation, the wind farm injection power is increased by 0.1 p.u. at 2 s and by 0.3 p.u. at 3 s. The time-domain simulation results for the original system, the system with the traditional damping controller, and the system with the DRDOC are presented in Fig. 10. When the first disturbance occurs, the system automatically recovers stability. However, after the second disturbance, due to its larger magnitude, the original system loses stability and begins to oscillate. The traditional damping controller improves stability to some extent, but small oscillations remain. In contrast, the system with the DRDOC quickly regains stability, highlighting the superior effectiveness of the DRDOC.

(a)

(b)

Fig. 10 Simulation with disturbances. (a) active power (b) DC voltage.

*D. Distributional robustness of DRDOC.*

This section evaluates the robustness of the DRDOC against the true distribution of disturbances. Based on the distribution in [24], 60% of the random disturbance points are selected within the ±0.3 p.u. range, while the remaining 40% are distributed within the ±1.0 p.u. range. This simulates various scenarios that approximate a normal distribution while retaining randomness. This setup simulates various scenarios that approximate a normal distribution while preserving randomness. By including both large and small disturbances, the setup effectively verifies the performance of the DRDOC. A total of 500 random disturbance scenarios are generated.

The solution of DRO satisfies the requirements with a confidence level of 1-*β*, where *β* is set to 3%. A comparison of the damping performance under DRDOC, SO, and RO across these 500 scenarios is presented. In this paper, SO designs damping control parameters based on empirical distribution of disturbances. RO assumes the probability distribution of disturbances is unknown and describes random disturbances using ranges, optimizing control parameters under the worst-case disturbance scenario. The damping performance based on DRDOC is shown in Fig. 11.

In the worst-case scenario defined by the ambiguity set, DRDOC is required to achieve the damping probability exceeding 97%. Among these scenarios, DRDOC successfully ensures damping greater than zero in 97.4% of cases, surpassing the 97% target. This improvement is since the PDF of power fluctuations in the generated scenarios is more favorable than the worst-case distribution within the ambiguity set. However, a small fraction of extreme samples lies outside the ambiguity set's requirements, resulting in 2.6% of cases where the damping value falls below zero.

Fig. 11 System damping in 500 scenarios, based on DRDOC.

In the SO-based strategy, a specific probability distribution is assigned to the random disturbances, as described in [9], and SO is solved to determine the optimized parameters. Using the parameters obtained through SO, the damping performance is presented in Fig. 12. Among the 500 scenarios, the probability of achieving positive damping is 94.6%. However, SO is highly dependent on the assumed disturbance probability distribution. Since actual power fluctuations are often unpredictable, this leads to weaker performance under uncertain conditions.

RO requires a predefined disturbance range to solve the min-max optimization problem. In the RO-based strategy, the probability of achieving positive damping across the 500 scenarios is 93.0%, as shown in Fig. 13. The performance of RO is significantly influenced by the predefined disturbance range. While a larger disturbance range improves performance, it also increases control costs. Additionally, different disturbance ranges are configured for RO, with the corresponding results summarized in Table IV. As the disturbance range increases, the controller gain becomes larger, resulting in poorer economic performance in practical applications. Moreover, when the constraints on control parameters cannot be satisfied, RO fails to yield a solution.

In contrast, DRDOC incorporates the entire disturbance range as a probability distribution, achieving acceptable control costs while meeting damping requirements. When random disturbances follow a probability distribution, the control performance of the RO-based strategy is inferior to that of SO and DRDOC. Furthermore, as DRDOC offers distributional robustness compared to SO, it demonstrates superior damping performance under uncertain true distributions of disturbances.

In practical engineering, the distribution of wind power exhibits probabilistic characteristics, but its true distribution often remains unknown. Consequently, DRDOC is better suited for such operating conditions, effectively enhancing the system's distributional robustness.

Fig. 12 System damping in 500 scenarios, based on the SO.

Fig. 13 System damping in 500 scenarios, based on RO.

TABLE IV
DISTURBANCES RANGES AND COSTS OF RO

| Range | Parameter $K_m$ | $\xi > 0$ in 500 scenarios (%) |
|---|---|---|
| [-0.15, 0.15] | 19.8 | 90.5 |
| [-0.2, 0.2] | 26.7 | 93.0 |
| [-0.25, 0.25] | Infeasible | / |
| [-0.3, 0.3] | Infeasible | / |

VI. CONCLUSION

This paper proposes a distributionally robust damping optimization control framework to address instability caused by random disturbances in power systems. The DRDOC framework integrates Sobol sensitivity analysis and PF for optimal controller placement and key parameter selection. It employs PCE to quantify the nonlinear relationship between damping and disturbances. Additionally, ambiguity sets are utilized to represent the true uncertainty distribution of disturbances. A simplified form of the DRO problem is derived to enhance computational feasibility. The DRDOC is versatile, with no specific limitations on system topology or operating modes, making it suitable for widespread application in multi-terminal connection scenarios. Depending on the actual disturbance conditions, DRDOC can be deployed on different converters, offering greater flexibility. Compared to traditional damping control, DRDOC effectively mitigates oscillations caused by a variety of random disturbances. Furthermore, DRDOC outperforms both RO and SO by addressing scenarios where disturbances exhibit probabilistic characteristics, but their true distributions are unknown. This ensures superior damping performance and significantly enhances distributional robustness.

## Appendix A Base value of the VSC-MTDC

TABLE I
PARAMETERS OF THE VSC-MTDC SYSTEM

| Parameter name (unit) | Parameter value |
| --- | --- |
| VSC rated voltage (kV) | $\pm 200$ |
| VSC rated power (MW) | 450 |
| AC rated voltage (kV) | 220 |
| Converter loss (p.u.) | 1.5% |
| Transformer reactance (p.u.) | 0.15 |
| Transformer X/R ratio | 30 |
| capacities of WFVSC (MW) | 300 |
| Transmission line length of line1,line2, line3(km) | 50 100 50 |
| DC transmission line Resistance ($\Omega$/km) | 0.0113 |
| DC transmission line Inductance (mH/km) | 0.45 |
| DC transmission line Capacitance (0.28μF/km) | 0.28 |